\documentclass[11pt,reqno]{amsart}
\usepackage[foot]{amsaddr}

\title{Structure of lower tails in sparse random graphs}

\author{Byron Chin}
\address{Department of Mathematics, Massachusetts Institute of Technology}
\email{byronc@mit.edu}

\newcommand{\LT}{\mathcal{L}}
\newcommand{\q}{\mathbf{q}}
\numberwithin{equation}{section}

\usepackage{arxiv}

\begin{document}

\begin{abstract}%
    We study the typical structure of a sparse Erd\H{o}s--R\'enyi random graph conditioned on the lower tail subgraph count event. We show that in certain regimes, a typical graph sampled from the conditional distribution resembles the entropy minimizer of the mean field approximation in the sense of both subgraph counts and cut norm. The main ingredients are an adaptation of an entropy increment scheme of Kozma and Samotij, and a new stability for the solution of the associated entropy variational problem. 
    
    The proof can be interpreted as a structural application of the new probabilistic hypergraph container lemma for sparser than average sets, and suggests a more general framework for establishing such typical behavior statements.
\end{abstract}

\maketitle

\section{Introduction}
We study the Erd\H{o}s--R\'enyi random graph conditioned on the large deviation event of a lower tail subgraph count. For a fixed graph $H$, let $N_H(G)$ be the number of copies of $H$ (as a subgraph) in a graph $G$. We write $N_H(p)$ to denote the random variable $N_H(G(n,p))$, where $G$ is drawn from the Erd\H{o}s--R\'enyi distribution. The key question is: Given that the lower tail event $\LT_p(H, \eta) \vcentcolon= \{G: N_H(G) \leq \eta \E{N_H(p)}\}$ holds for some $\eta < 1$, what does a typical sample from the conditional distribution look like? 

The study of this question was initiated by the seminal work of Chatterjee and Varadhan \cite{CV:11}, who proved a general large deviations principle for the Erd\H{o}s--R\'enyi random graph with constant density. Their results imply the following two estimates written below, which we first describe informally. The first statement (\ref{eq: CV1}) says that the leading order in the exponent of the lower tail probability is determined by the way to induce few copies of $H$ with minimum entropic cost. The second statement (\ref{eq: CV2}) gives a more detailed description of the structure of the event, and says that the probability that a random sample with few copies of $H$ looks very different from these entropy minimizers is exponentially small. They showed
\begin{gather}
    -\log \Prob{\LT_p(H, \eta)} = (1+o(1))\Phi_p(H, \eta) \label{eq: CV1} \\
    \Prob{\delta_\square(G, M) > \epsilon ~ \vert ~ \LT_p(H, \eta)} \leq e^{-cn^2}. \label{eq: CV2}
\end{gather}
where $\Phi_p(H, \eta)$ is a variational problem that minimizes the entropic cost over the lower tail event (see Definition \ref{def: variational problem} for the precise formulation), $M$ is the set of minimizers to the variational problem $\Phi_p(H, \eta)$ in the space of graphons, and $\delta_\square$ is the cut distance. See \cite{CV:11} for a formal discussion of these expressions, or \cref{section: preliminaries} for the definitions and their role in this paper. Much finer questions have also been studied in the case of dense graphs, determining the typical structure of large graphs with various fixed subgraph densities (for example \cite{NRS:23, KRRS:18, Rad:18} and references within).  

There are two barriers that prevent the extension of \cite{CV:11} to the case of sparse random graphs. First, their proof relies on Szemer\'edi's regularity lemma which is ineffective at handling graphs with vanishing edge density. Second, and perhaps more significantly, the large deviation principle is formulated in terms of graphons, a natural limit object for dense graphs. The space of graphons is a compact topological space which nicely encodes similarity and distance between graphs, see \cite{Lov:12} for a comprehensive survey. Theories for the limit of sequences of sparse graphs exist, e.g. \cite{BCCZ:19}, but there is no convenient limit object which captures all the same properties as a graphon. Consequently it is not clear how a corresponding large deviation principle should be formulated. Despite these difficulties with generalizing the entire large deviations principle, it is possible to study these tail probabilities in a more specific manner. 

For equation (\ref{eq: CV1}), a natural first step is to replace Szemer\'edi's regularity lemma with the weak regularity lemma of Kannan and Frieze \cite{FK:99} which can handle slowly vanishing densities $p \geq (\log n)^{-c}$. A breakthrough technique for computing large deviation probabilities for nonlinear functions of independent Bernoulli random variables was introduced by Chatterjee and Dembo \cite{CD:16}, which extends the range of densities to $p \geq n^{-\alpha(H)}$ for some explicit function $\alpha(H)$. This inspired a line of work developing the technique for arbitrary graphs and even uniform hypergraphs \cite{Eld:18, CD:20, CDP:21}. Recently, a work of Kozma and Samotij \cite{KS:23} showed (\ref{eq: CV1}) for $p \gg n^{-1/m_2(H)}$, the full range of densities for which this statement can be expected to hold, using an entropy increment technique. Here $m_2(H)$ refers to the 2-density of a graph $H$ defined in \cref{section: preliminaries}. For sparser graphs, the lower tail probability is determined by a Poisson-type behavior, and Janson's inequality gives a tight bound. Janson and Warnke \cite{JW:16} showed that this is true in a more general setting -- see the references within for an overview of work in this even sparser regime. These lines of work present a fairly complete picture of the lower tail probability. 

Extending the structural statement (\ref{eq: CV2}) to the sparse setting has seen far less progress. One reason for this is that the variational problem for lower tail subgraph counts has proven to be difficult to analyze. Very little is known about the set of minimizers, so sparse analogs of (\ref{eq: CV2}) are more difficult to show explicitly. Zhao showed that for $\eta$ sufficiently close to 1 the constant function is the unique minimizer of $\Phi_p(H, \eta)$, whereas for $\eta$ sufficiently close to 0 the constant function is no longer a minimizer \cite{Zha:17}. He studied a sparse limit of the variational problem, leading his results to hold for $p \to 0$. This remains the best known range of parameters for which a solution to the lower tail variational problem is known. To the best of our knowledge, no analogs of (\ref{eq: CV2}) have been shown for the lower tail for any $p = o(1)$. 

\subsection{Main Results}
We provide a framework for deducing sparse analogs to (\ref{eq: CV2}) -- the typical structure of a conditioned random graph -- from a solution to the variational problem. We apply this argument in the setting of lower tail subgraph counts of $G(n,p)$, where all the ingredients are available to prove the complete result. We show that conditioned on a random graph having few copies of one graph $H'$, subgraph counts of another graph $H$ are accurately predicted by the solution to the variational problem with high probability. 
\begin{theorem}\label{thm: main intro}
    Fix any graphs $H'$ and $H$. Let $G_\LT \sim G(n, p)$ conditioned on $\LT_p(H', \eta)$, where $\eta > \eta_{H'}$. Let $q = \eta^{1/e(H')}p$ and $m = \max\{m_2(H), m_2(H')\}$. For any $\epsilon > 0$ there exists an $L(\epsilon, \eta, H, H')$ such that if $p \geq Ln^{-1/m}$, 
    \[ \Prob{|N_{H}(G_\LT) - \E{N_{H}(q)}| > \epsilon \E{N_{H}(q)}} < \epsilon. \]
\end{theorem}
\begin{remark}
    For $p \gg n^{-1/m}$ carefully tracking asymptotics through the argument of \cite{KS:23} should yield a decaying bound. However, as written the proof only obtains a bound for arbitrary fixed $\epsilon$. 
\end{remark}
Here $\eta_{H'}$ is the constant from \cite{Zha:17} above which the solution to the variational problem is known to be constant. Numerically $\eta_{H'}$ is defined as the solution to a fixed point equation and can be approximated. For example, $\eta_{K_3}$ is $0.1012...$ -- see \cite[Section 5]{Zha:17} for more estimates.

As a consequence, we deduce typical structure of the graph in the sense of cut norm as well. The cut norm is a natural notion of distance between graphs that enjoys a set of nice properties when extended to the space of graphons (see \cref{def: cut norm}).
\begin{corollary}\label{cor: main cor}
    Fix any graph $H$. Let $G_\LT \sim G(n, p)$ conditioned on $\LT_p(H, \eta)$ where $\eta > \eta_H$. Let $q = \eta^{1/e(H)}p$. For any $\epsilon > 0$ there exists an $L(\epsilon, \eta, H)$ such that if $p \geq Ln^{-1/m_2(H)}$, 
    \[ \Prob{\norm{G_\LT - q}_\square > \epsilon p} < \epsilon. \]
\end{corollary}

Our results extend (\ref{eq: CV2}) to the full range of densities for which the variational problem is relevant, and moreover establish typical structure in the stronger sense of subgraph counts. In the sparse setting, it is known that subgraph counts control the cut norm, but the converse fails. Given our results about the typical structure of $G(n, p)$ conditioned on a lower tail event, one might wonder if an even stronger structural relationship holds. One possible formulation is the following:
\begin{question}
    Fix any graph $H$, $\eta > \eta_H$ and let $q = \eta^{1/e(H)}p$. Let $p \gg n^{-1/m_2(H)}$ and $G_{\LT} \sim G(n, p)$ conditioned on $\LT_p(H, \eta)$. For any $\epsilon > 0$, does there exist a coupling such that with high probability
    \[ G(n, (1-\epsilon)q) \subset G_\LT \subset G(n, (1+\epsilon)q)? \]
\end{question}
A first step to investigate the above question is to study various statistics of the lower tail event that would follow from the existence of such a coupling.
\begin{question}
    Let $G_{\LT} \sim G(n, p)$ conditioned on $\LT_p(H, \eta)$. What is the size of the largest clique in $G_{\LT}$?
\end{question}

\subsection{Proof Overview}
We combine a probabilistic version of the hypergraph container lemma with a new stability result to deduce the typical structure.

Kozma and Samotij \cite{KS:23} introduced a new entropy-based approach to estimating lower tail probabilities. They interpret their work as a ``weak, probabilistic version of the hypergraph container lemma for sparser-than-average sets". The proof relies on an entropy increment iteration to obtain approximate independence within a random graph conditioned on having few copies of $H'$. In particular, by conditioning on a small subset of edges, they decompose the lower tail measure into a mixture of measures which behave nicely with respect to copies of $H'$. We leverage their iteration scheme to decompose the measure with respect to multiple graphs, and obtain approximate independence with respect to each of these graphs simultaneously.

We extend the argument of Zhao \cite{Zha:17} for showing that the constant function is the unique entropy minimizer when the lower tail parameter $\eta$ is sufficiently close to 1. We obtain a corresponding stability result stating that near minimizers must look like the constant function.

Our proof then proceeds as follows. Given the distribution of $G(n, p)$ conditioned on $\LT_p(H', \eta)$, condition on a small set of edges to obtain control over copies of both $H$ and $H'$. Then, with high probability, an independent sample from the marginal probabilities of each edge from this conditioned distribution will also satisfy a lower tail event for $H'$. We then invoke our new stability result to deduce that the marginal edge probabilities must be close to the minimizer of the variational problem. From this we deduce that the expected number of copies of $H$ in this independent sample is close to the quantity predicted by the minimizer. Recalling our control on copies of $H$ as well, we deduce that the expected number of copies of $H$ under the joint distribution of edges is close to the independent sample from the marginals, and thus is also accurately predicted by the solution to the variational problem. Finally, to boost this estimate on the expectation to obtain concentration, we leverage the fact that the variance can be expressed in terms of subgraph counts, and apply a second moment argument.

The bulk of this argument is in fact very general, and easily allows for control over multiple subgraph counts as well as conditioning on multiple lower tail events. Furthermore, the main ideas can be phrased in the language of counting edges in hypergraphs, a flexible viewpoint exemplified in the transference principles of Conlon and Gowers \cite{CG:16} and Schacht \cite{Sch:16} and the hypergraph container theorems \cite{BMS:15, ST:15}. This would allow one to deduce a ``typicality" statement for edge counts of two hypergraphs living on the same vertex set, provided we have the necessary stability of the variational problem. The primary bottleneck to these types of results appears to be the analysis of the entropy variational problem. 

\subsection{Organization} 
The next section sets up the definitions and notation for the remaining arguments. In \cref{section: entropy increment} we recall the main results of the entropy increment scheme to obtain approximate independence. In \cref{section: stability} we prove the new stability of the entropy variational problem. \cref{section: main proof} combines the ingredients to complete the proof of the main theorem. We conclude in Section \ref{section: discussion} discussing our argument in a more general setup.

\subsection*{Acknowledgments}
The author thanks Yufei Zhao for introducing the problem and many conversations about this work. The author also thanks Wojciech Samotij and Mehtaab Sawhney for helpful discussions. The author is supported by an NSF Graduate Research Fellowship.

\section{Preliminaries}\label{section: preliminaries}
In this section we introduce the primary definitions and notation that will be used throughout the paper. Some notation will be imported from \cite{KS:23} and \cite{Zha:17}, and we will give a reference when this is the case. From here until Section \ref{section: discussion} we will restrict our attention to the setting of subgraphs of random graphs. The more general setup will be introduced in the discussion in Section \ref{section: discussion}. 

We let $G = (V, E)$ denote a graph and record the following statistics about fixed graphs.
\begin{definition}[Graphs]\label{def: graphs}
    Let $v(G) = |V|$ be the number of vertices and $e(G) = |E|$ be the number of edges of of $G$. We use $N_H(G)$ to denote the number of copes of $H$ as a subgraph of $G$. A key statistic is the 2-density of $G$, which is defined as the following
    \[ m_2(G) \vcentcolon= \max \left\{ \frac{e(F) - 1}{v(F) - 2}: F \subset G, e(F) \geq 2\right\}. \]
\end{definition}
In this paper we are interested in the behavior of random graphs, which will be notated according to the  definition below. 
\begin{definition}[Random graphs]\label{def: random graphs}
    For a vector $\q \in [0,1]^{\binom{n}{2}}$ let $G(n, \q)$ be the random graph with vertex set $[n]$ where each edge $e = (i,j)$ is included independently with probability $q_{i,j} = q_e$. When $\q = p$ is a constant vector, we simply write $G(n,p)$, the well-known Erd\H{o}s--R\'enyi random graph. We use $N_H(G)$ to represent the number of copies of $H$ as a subgraph of $G$. We use $N_H(\q)$ to represent $N_H(G(n, \q))$, the (random) number of copies of $H$ in a sample from $G(n, \q)$. 
\end{definition}
We are primarily interested in the lower tail subgraph count event, when the random graph has much fewer copies of a subgraph $H$ than expected. 
\begin{definition}[Lower tail event]\label{def: lower tail}
    Let $\LT_p(H, \eta)$ denote the event that $N_H(G) \leq \eta\E{N_H(p)}$. We write $G_\LT$ for a graph drawn from $G(n, p)$ conditioned on $\LT_p(H, \eta)$.
\end{definition}
One notion of similarity we use to characterize the distance between instances of the random graphs is the cut norm. For our purposes the cut norm can be defined on the set of $n \times n$ matrices as follows:
\begin{definition}[Cut norm]\label{def: cut norm}
    For $A \in \mathbb{R}^{n \times n}$,
    \[ \norm{A}_\square \vcentcolon= \sup_{x, y \in [0,1]^n} \frac{1}{n^2}\abs{x^\top A y}. \]
    We also use the cut norm on vectors $\q \in [0,1]^{\binom{n}{2}}$. In this case we interpret $\norm{\q}_\square = \norm{A(\q)}_\square$ where $A(\q)_{ij} = A(\q)_{ji} = q_{i,j}$. 
\end{definition}
One should think of $A$ as the difference of adjacency matrices of two (weighted) graphs for this paper. In this case, we obtain the cut norm on graphs alluded to earlier. For a detailed account of the full definition and properties of the cut norm in relation to the theory of graph limits, see the monograph by Lov\'asz \cite{Lov:12}. 

Many of the ideas in the paper can also be conveniently phrased in the language of hypergraphs, as done in \cite{KS:23}. A hypergraph will be represented by $\mathcal{H} = (\mathcal{V}, \mathcal{E})$. We similarly use $v(\mathcal{H}) = |\mathcal{V}|$ and $e(\mathcal{H}) = |\mathcal{E}|$. 
\begin{definition}[Hypergraph of copies of $H$]\label{def: hypergraph}
    Let $\mathcal{H}(H)$ be the hypergraph with vertex set $\mathcal{V} = \binom{[n]}{2}$ and edge set $\mathcal{E} = $ copies of $H$ as a subgraph of the complete graph $K_n$. We will use $\mathcal{H}(H)$ to refer to both the hypergraph and its corresponding edge set.
\end{definition}
The following vector representation of a vertex set will be convenient for us to work with. 
\begin{definition}[Vector representation]\label{def: vector}
    Let $R_p$ be a random subset of $\mathcal{V}$ where each vertex is included independently with probability $p$. Let $Y \in \{0,1\}^{\binom{[n]}{2}}$ be the indicator vector of this random subset conditioned on $\LT_p(H, \eta)$, so that $Y$ corresponds to the edges of $G_\LT$. We write $Y_A = \prod_{e \in A} Y_e$ for the indicator that the subset $A$ is included in $R_p$. 
\end{definition}

The relative entropy between random variables plays a key role, as it controls the probability of the lower tail event. Here we define the relevant quantities and list a few useful properties. See \cite[Section 4]{KS:23} for a more in depth introduction. 
\begin{definition}[Relative entropy]\label{def: entropy}
    For a binary vector $X$ define its $p$-divergence by 
    \[ I_p(X) = D_\mathrm{KL}(X \ || \ \mathrm{Ber}(p)^k).\] 
    We also define the conditional divergence by $I_p(X|Z) = \E{I_p(X^Z)}$ where $X^Z$ denotes the random variable $X$ conditioned on $Z$. For any binary random vectors $X_1, X_2$ and any random variable $Z$, we have the following useful properties:
    \begin{enumerate}[label=\arabic*.]
        \item $I_p(X_1) \geq 0$.
        \item $I_p(X_1,X_2) = I_p(X_1|X_2) + I_p(X_2) \geq I_p(X_1) + I_p(X_2)$.
        \item $I_p(X_1, X_2|Z) \geq I_p(X_1|Z) + I_p(X_2|Z)$. 
        \item $I_p(Y) = -\log \Prob{\LT_p(H, \eta)}$ where $Y$ is the conditional vector from Definition \ref{def: vector}.
    \end{enumerate}
    
    In the case that $X$ itself is a Bernoulli variable we recover the relative entropy function. Let $i_p(q) = q \log \frac{q}{p} + (1-q)\log\frac{1-q}{1-p}$ be the relative entropy of $q$ with respect to $p$. For a vector $\q$, we will write $i_p(\q) \vcentcolon= \sum_{e \in \binom{[n]}{2}} i_p(q_e)$ for the total relative entropy with respect to $p$. 
\end{definition}
The following variational problem encodes the relationship between entropy and the lower tail event.
\begin{definition}[Variational problem]\label{def: variational problem}
    We define the variational problem associated with the lower tail event $\LT_p(H, \eta)$ as
    \[ \Phi_p(H, \eta) \vcentcolon= \min_{\q \in [0,1]^{\binom{n}{2}}}\left\lbrace i_p(\q): \E{N_H(\q)} \leq \eta \E{N_H(p)}\right\rbrace. \]
\end{definition}
We follow the notation used by Kozma and Samotij \cite{KS:23}, but note here that the relative entropy function is denoted by $I_p(q)$ rather than $i_p(q)$ and the variational problem is denoted by $\mathrm{LT}_p(H, \eta)$ rather than $\Phi_p(H, \eta)$ in \cite{Zha:17}.

\section{Entropy Increment}\label{section: entropy increment}
In this section we collect a series of results from the work of Kozma and Samotij \cite{KS:23} that controls the probability of $\LT_p(H, \eta)$ via an entropy increment argument alluded to in the introduction.

The first lemma says that conditioned on a lower tail event, the marginal probabilities that fixed sets of edges appear can only decrease. The proof is an application of Harris' inequality.
\begin{lemma}[{\cite[Claim 16]{KS:23}}]\label{lemma: harris}
    For every $W \subset \mathcal{V}$ and $A \subset \mathcal{V} \setminus W$, $\E{Y_A \vert (Y_w)_{w \in W}} \leq p^{|A|}$.
\end{lemma}

A key step in \cite{KS:23} shows approximate conditional independence of copies of $H$ under the lower tail event by conditioning on a small subset of the edges. This translates to a conditional independence of the appearance of edges of $\mathcal{H}(H)$ given the configuration on a small subset of $\mathcal{V}$. The conditional independence is measured by the following quantities.  
\[ D_W(B, b) \vcentcolon= \E{Y_B \vert (Y_w)_{w \in W}} - \E{Y_{B \setminus \{b\}} \vert (Y_w)_{w \in W}}\E{Y_b \vert (Y_w)_{w \in W}} \]
measures the correlation of a single vertex $b \in \mathcal{V}$ with a hyperedge $B \in \mathcal{E}$ conditioned on the configuration of $W \subset \mathcal{V}$. The weighted sum of these squared correlations is captured by the following quantity:
\begin{equation}\label{eq: energy}
    \mathcal{E}_H(W) \vcentcolon= \E{\sum_{A \in \mathcal{H}(H) - W}\sum_{B \subset A, |B| \geq 2}\sum_{b \in B}\frac{D_W(B, b)^2}{\binom{r}{|B|}|B|p^{2|B|}}}.
\end{equation}  
Following \cite{KS:23}, we use $\mathcal{H}(H) - W$ denote the induced sub-hypergraph of $\mathcal{H}(H)$ on the vertex set $\mathcal{V} \setminus W$. It turns out that this quantity controls the difference between a sample from the joint distribution of $Y$ and sampling each edge independently with the corresponding marginal. The proof is an application of inclusion-exclusion and the Cauchy--Schwarz inequality, see \cite[p.686]{KS:23} for the proof.
\begin{lemma}\label{lemma: cauchy schwarz}
    Let $e(H) = r$ and $W \subset \mathcal{V}$ be a subset of the vertices of $\mathcal{H}(H)$. 
    \[ \E{\sum_{A \in \mathcal{H}(H)-W} \abs{\E{Y_A\vert (Y_w)_{w \in W}} - \prod_{a \in A}\E{Y_a \vert (Y_w)_{w \in W}}}} \leq p^r((r-1)e(\mathcal{H}(H)-W))^{1/2}\mathcal{E}_H(W)^{1/2}. \]
\end{lemma}
Thus, to well-approximate the distribution of (hyper)edges in $Y$ from its independent marginals, we need a $W$ such that both $|W|$  and $\mathcal{E}_H(W)$ are small. The following generalization of \cite[Lemma 18]{KS:23} finds this small set, and follows from the same proof given in \cite{KS:23} and the discussion in \cite[Section 7]{KS:23}.

\begin{lemma}\label{lemma: increment}
    Let $H$ be a non-empty graph and $\mathcal{H} = \mathcal{H}(H)$. For all positive $\alpha$ and $\beta$, there exists $L$ and $V_0$ such that the following holds: if $|\mathcal{V}| \geq V_0$ and $p \geq Ln^{-1/m_2(H)}$ and $Y$ is a distribution on subsets of $\mathcal{V}$ such that $\E{Y_A \vert (Y_w)_{w \in W}} \leq p^{|A|}$ for all $A \subset \mathcal{V} \setminus W$, there exists a set $W$ with at most $\alpha v(\mathcal{H})$ elements such that $\mathcal{E}_H(W) \leq \beta e(\mathcal{H})$. 
\end{lemma}

Lemmas \ref{lemma: harris}, \ref{lemma: cauchy schwarz}, and \ref{lemma: increment} can be combined to obtain the following theorem, showing that the probability of the lower tail event is controlled to the leading term in the exponent by the associated variational problem. 
\begin{theorem}[{\cite[Theorem 2]{KS:23}}]\label{theorem: tail probability}
    For every non-empty graph $H$, $p_0 < 1$, and every $\epsilon > 0$, there exists a constant $L$ such that for $Ln^{-1/m_2(H)} \leq p \leq p_0$ and every $\eta \in [0,1]$, 
    \[ (1-\epsilon)\Phi_p(H, \eta+\epsilon) \leq -\log \Prob{\LT_p(H, \eta)} \leq (1+\epsilon)\Phi_p(H, \eta(1-\epsilon)). \]
\end{theorem}

\section{Stability of the variational problem}\label{section: stability}
In the previous section we related the lower tail probability to an entropy variational problem. The goal of this section is to give a quantitative description to the solution of the variational problem. Our main new technical ingredient is a refinement of  \cite{Zha:17} that enhances the uniqueness of the solution to $\Phi_p(H, \eta)$ to a stability of the minimizer. In this section, we work with the sparse limit of the variational problem, which is an appropriate substitute for $\Phi_p(H, \eta)$ when $n$ is sufficiently large and $p = o(1)$. We define $h(x) \vcentcolon= x\log x - x + 1$ so that $\lim_{p \to 0} p^{-1}i_p(px) = h(x)$ uniformly for $x \in [0,1]$. It follows that 
\[ \lim_{p \to 0} p^{-1}\Phi_p(H, \eta) = \Phi(H, \eta) \]
where 
\[ \Phi(H, \eta) \vcentcolon= \min_{\q \in [0,1]^{\binom{n}{2}}}\left\lbrace \sum_e h(q_e): \E{N_H(\q)} \leq \eta \cdot N_H(1) \right\rbrace. \]
We know from \cite{Zha:17} the following theorem:
\begin{theorem}\label{thm: minimizer}
    For every graph $H$ there exists an $\eta_H < 1$ such that for all $\eta > \eta_H$, $\Phi(H, \eta)$ is uniquely minimized by the constant function, in this case $\q \equiv \eta^{1/e(H)}$. 
\end{theorem}
We show that the constant function is stable as the minimizer to $\Phi(H, \eta)$ in this regime in the following quantitative sense.
\begin{proposition}\label{prop: stability}
    Let $H$ be any graph. There exists an $\eta_H$ such that for all $\eta > \eta_H$ there exists a $C_{\eta, H}$ such that the following holds. Let $q = \eta^{1/e(H)}$ and $\epsilon > 0$ be sufficiently small in terms of $\eta$ and $H$. Then $\q \equiv q$ is the unique minimizer of $\Phi(H, \eta)$, and if $\sum_e h(q_e) < \Phi(H, \eta) + \epsilon^{12}\binom{n}{2}$, then $\norm{\q - q}_\square < C_{\eta, H}\epsilon$.
\end{proposition}
To begin, we first give a weak description of a near minimizer, showing that only few of its entries can be very small. We will use subgraph count density notation in this section, writing $t_{\mathrm{inj}}(H, \q) = \frac{\E{N_H(\q)}}{N_H(1)}$. This is also known as the injective homomorphism density, closely related to the standard notion of homomorphism density (see \cite[Section 4.3]{Zha:23} for more on homomorphism densities). The following is the stability version of \cite[Lemma 5.2]{Zha:17}.

\begin{lemma}\label{lemma: stability 2}
    Let $\eta > \eta_H$. For all $\epsilon$ sufficiently small in terms of $\eta, H$, if $\sum_e h(q_e) \leq \binom{n}{2}h(\eta^{1/e(H)}) + \epsilon^2\binom{n}{2}$ and $t_\mathrm{inj}(H, \q) < \eta$ then $q_e \geq c \vcentcolon= \eta^{1/\eta}$ for all but $\epsilon\binom{n}{2}$ edges $e$. 
\end{lemma}
\begin{proof}
    The idea is to show that if the set of small entries is too large in a near minimizer, then a marginal boost of the edge probabilities on this set leads to a configuration with few copies of $H$ that has too small of an entropy.
    
    Let $B = \{e: q_e < c\}$ denote the bad set of edges and let $\delta = \frac{|B|}{\binom{n}{2}}$. Suppose we increment $\q$ on each of these edges to obtain $\q' \vcentcolon= \q + \gamma \mathbbm{1}_B$. Since $h$ is convex and decreasing, for every $e \in B$ we have
    \[ h(q_e + \gamma) - h(q_e) \leq h(c + \gamma) - h(c) = -\gamma(1 - \log(c+\gamma)) - c\log(1 - \frac{\gamma}{c+\gamma}) < 0. \]
    Setting $q \vcentcolon= \eta^{1/e(H)}$, we can bound the new entropy by 
    \begin{align*}
        \sum_e h(q_e') &= \sum_{e \in B} h(q_e + \gamma) + \sum_{e \in B^c} h(q_e) \\
        &\leq \sum_e h(q_e) + |B|\left(-\gamma(1 - \log(c+\gamma)) - c\log(1 - \frac{\gamma}{c+\gamma})\right) \\
        &\leq \binom{n}{2}\left(h(q) + \epsilon^2 -\delta\gamma(1 - \log(c+\gamma)) - \delta c\log(1 - \frac{\gamma}{c+\gamma}) \right)
    \end{align*}
    On the other hand, since $\norm{\q' - \q}_\square = \norm{\gamma\mathbbm{1}_B}_\square < \gamma\delta$, the counting lemma (see e.g. \cite[Section 4.5]{Zha:23}) implies 
    \[ t_{\mathrm{inj}}(H, \q') - t_{\mathrm{inj}}(H, \q) \leq e(H)\norm{\q' - \q}_\square \leq \gamma\delta e(H). \]
    In particular, 
    \[ t_{\mathrm{inj}}(H, \q') \leq \eta + \gamma\delta e(H). \]
    By \cref{thm: minimizer} we know that the variational problem $\Phi(H, \eta + \gamma\delta e(H))$ is minimized by the constant function $(\eta + \gamma\delta e(H))^{1/e(H)}$, implying that
    \begin{align*}
        \sum_e h(q_e') &\geq \binom{n}{2} h((\eta + \gamma\delta e(H))^{1/e(H)}) \geq \binom{n}{2}h(q + \gamma \delta \eta^{-1}) \\
        &= \binom{n}{2}\left(h(q) - \frac{\gamma\delta}{\eta}\left(1 - \log(q+\frac{\gamma\delta}{\eta})\right) - q\log(1 - \frac{\gamma\delta \eta^{-1}}{q+\gamma\delta \eta^{-1}})\right).
    \end{align*}
    Combining the two bounds on $\sum_e h(q_e')$, we must have 
    \begin{align*}
        - \frac{\gamma\delta}{\eta}\left(1 - \log(q+\frac{\gamma\delta}{\eta})\right) - q\log(1 - \frac{\gamma\delta \eta^{-1}}{q+\gamma\delta \eta^{-1}}) \leq \epsilon^2 -\delta\gamma(1 - \log(c+\gamma)) - \delta c\log(1 - \frac{\gamma}{c+\gamma})
    \end{align*}
    By Taylor expansion, we know 
    \[ \delta c\log(1 - \frac{\gamma}{c + \gamma}) = -\frac{\delta c\gamma}{c + \gamma} - O_\eta(\delta \gamma^2) = -\gamma\delta - O_\eta(\delta\gamma^2) \]
    and similarly 
    \[ q\log(1 - \frac{\gamma\delta \eta^{-1}}{q+\gamma\delta \eta^{-1}}) = - \frac{\gamma\delta\eta^{-1}}{1 + \gamma\delta (q\eta)^{-1}} - O_\eta(\delta^2\gamma^2) = -\gamma\delta \eta^{-1} - O_\eta(\delta^2\gamma^2). \]
    Rearranging, we obtain 
    \begin{gather*}
        \gamma\delta\left( 1 - \log(c + \gamma) - \frac{1}{\eta} + \frac{1}{\eta}\log(q + \frac{\gamma\delta}{\eta}) - 1 + \frac{1}{\eta} \right) + O_\eta(\delta \gamma^2) \leq \epsilon^2 \\
        \gamma\delta\left( - \log c + \frac{1}{\eta}\log q \right) + O_\eta(\delta \gamma^2) \leq \epsilon^2
    \end{gather*}
    Setting $\gamma = \frac{2\epsilon}{\frac{\log q}{\eta} - \log c}$, we deduce that for all $\epsilon$ sufficiently small in terms of $\eta$ and $H$, the number of bad edges is  $\leq \epsilon\binom{n}{2}$. 
\end{proof}

Before proving the stability, we record a useful fact that appears as \cite[Fact 5.5]{Zha:17}. 
\begin{fact}\label{fact: stability}
    Let $\eta_H$ be the unique solution in the interval $(0,1)$ to the equation 
    \[ h(\eta^{1/\eta}) = h(\eta^{1/e(H)}) + \eta^{1/e(H)}\log(\eta^{1/e(H)})(\log(\eta^{1/\eta}) - \log(\eta^{1/e(H)})). \]
    Then for all $(x, r) \in [\eta_H^{1/\eta_H}, 1] \times [\eta_H^{1/e(H)}, 1]$ the inequality $h(x) \geq h(r) + r\log(r)(\log x - \log r)$ holds.
\end{fact}

\begin{proof}[Proof of \cref{prop: stability}]
    The strategy is as follows. We use \cref{lemma: stability 2} and \cref{fact: stability} to give a lower bound on the average value of $\log q_e$. Then, a nearly matching lower bound follows by Jensen's inequality. The combination of these two bounds implies that each inequality must be nearly tight, from which we can extract term-by-term bounds on $q_e$. 

    Applying \cref{lemma: stability 2} we know that $B = \{e: q_e < \eta^{1/\eta}\}$ has size at most $\epsilon^6 \binom{n}{2}$. We show that for most of the remaining edges, that $q_e$ must be close to $q \vcentcolon= \eta^{1/e(H)}$. Using \cref{fact: stability}, we lower bound 
    \begin{align}\label{eq: term-wise}
        \sum_e h(q_e) &\geq \sum_{e \in B^c} h(q_e) \geq \sum_{e \in B^c} \left( h(q) + q\log q(\log q_e -\log q)\right) \\
        &= |B^c|h(q) + q\log q\left(\sum_{e \in B^c}(\log q_e - \log q)\right)
    \end{align}
    Combining this with the assumption that $q_e$ is a near minimizer, we obtain
    \[ \sum_{e \in B^c} (\log q - \log q_e) \leq \frac{\epsilon^6 + \epsilon^{12}}{q \log \frac{1}{q}} \binom{n}{2} \] 
    so that 
    \[ \frac{1}{\binom{n}{2}}\sum_{e \in B^c} \log q_e \geq (1-\epsilon^6) \log q - \frac{\epsilon^6 + \epsilon^{12}}{q \log \frac{1}{q}} = \log q - O_{\eta, H}(\epsilon^6). \]

    We now work out the upper bound. Recall that $N_H(1)$ is the total number of copies of $H$ in $K_n$. By symmetry, each edge is in $\frac{e(H)}{\binom{n}{2}}N_H(1)$ copies of $H$. Thus, the number of copies of $H$ that intersect $B$ is at most $\epsilon^6 e(H)N_H(1)$. Let $X$ be a uniformly random copy of $H$ that is contained in $B^c$. Then
    \begin{align*}
        \E{\log(\prod_{e \in X} q_e)} &= \frac{1}{N_H(\mathbbm{1}_{B^c})}\sum_{H} \sum_{e \in H} \log q_e = \frac{1}{N_H(\mathbbm{1}_{B^c})} \sum_{e \in B^c} \log q_e \cdot \abs{\{H: e \in H\}} \\
        &\geq \frac{1}{N_H(\mathbbm{1}_{B^c})} \sum_{e \in B^c} \frac{e(H)}{\binom{n}{2}}N_H(1)\log q_e \geq \frac{e(H)}{\binom{n}{2}(1-\epsilon^6 e(H))} \sum_{e \in B^c} \log q_e.
    \end{align*}
    On the other hand, Jensen's inequality implies 
    \[ \E{\log(\prod_{e \in X} q_e)} \leq \log\E{\prod_{e \in X} q_e} \leq \log(\frac{\eta}{1-\epsilon^6 e(H)}) = e(H)\log q - \log(1 - \epsilon^6 e(H)). \]
    Combining these inequalities yields 
    \[ \frac{1}{\binom{n}{2}}\sum_{e \in B^c} \log q_e \leq (1-\epsilon^6 e(H))\log q - \frac{1-\epsilon^6 e(H)}{e(H)}\log(1-\epsilon^6 e(H)) = \log q + O_{\eta, H}(\epsilon^6). \]

    The lower and upper bounds imply that the term-wise application of \cref{fact: stability} in \cref{eq: term-wise} must be nearly tight. Quantitatively, the number of edges $e$ for which 
    \[ h(q_e) > h(q) + q\log q(\log q_e - \log q) + \epsilon^{3} \]
    is at most $O_{\eta, H}(\epsilon^{3})\binom{n}{2}$. One can check that the equation $h(q_e) = h(q) + q\log q(\log q_e - \log q)$ has at most two solutions on $[0,1]$, one of which is $q_e = q$. 

    We analyze the order of vanishing of the function $h(q_e)-h(q)-q\log q(\log q_e-\log q)$ at each of the roots. Taking two derivatives, we get the expression $\frac{1}{q_e} + \frac{q\log q}{q_e^2}$, which vanishes at $q_e = q\log \frac{1}{q}$. On the other hand, the third derivative is $-\frac{1}{q_e^2} - \frac{2q\log q}{q_e^3}$, which vanishes at $q_e = 2q\log \frac{1}{q}$. In particular, since $0 < q < 1$, one of these two derivatives does not vanish at every point $q_e$. This implies that the order of vanishing at each of the roots is at most 3.
    
    Thus, for $\epsilon$ sufficiently small in terms of $\eta$ and $H$, if $h(q_e) \leq h(q) + q\log q(\log q_e - \log q) + \epsilon^{3}$ we must have that $q_e$ is within $O(\epsilon)$ of one of these solutions. Finally, since we know that the average value of $\log q_e$ is within $O(\epsilon^6)$ of $\log q$, all but an $O(\epsilon^6)$ fraction of the $q_e$ must in fact be close to $q$ rather than the second root. This shows that at most $O(\epsilon^{3}\binom{n}{2})$ of the edges $q_e$ have $|q_e - q| > \epsilon$, which implies that $\norm{\q - q}_\square < C_{\eta, H}\epsilon$.
\end{proof}

\section{Proof of the main theorem}\label{section: main proof}
In this section we prove the main results Theorem \ref{thm: main intro} and Corollary \ref{cor: main cor}. For the rest of the section we retain the following setup: Fix a graph $H'$ and suppose $\eta > \eta_{H'}$. Let $G_{\LT} \sim G(n, p)$ conditioned on $\LT_p(H', \eta)$. Let $q = \eta^{1/e(H')}p$ be the unique minimizer of $\Phi_p(H', \eta)$. 

The following lemma establishes a weak notion of typical structure by showing that the conditional distribution can be decomposed into a mixture of a small number of measures, almost all of which have marginals that are close to the optimizer of the variational problem. 
\begin{lemma}\label{lemma: weak structure}
    For any $\epsilon > 0$ there exists an $L$ such that if $p \geq Ln^{-1/m_2(H')}$ then the following holds. There exists a set of edges $W_0$ such that for any $W \supseteq W_0$ and $|W| < \frac{\epsilon}{e(H')}\binom{n}{2}$ we have 
    \[ \Prob{\norm{\q^W - q}_\square > \epsilon p} \lesssim \epsilon \]
    where $q^W_e = \E{Y_e \vert (Y_w)_{w \in W}}$ for $e \not \in W$ and $q^W_e = p$ otherwise.
\end{lemma}
\begin{proof}
    Throughout this proof, we may assume that $\epsilon$ is sufficiently small as a function of $\eta$ and $H'$. Set $\alpha = \frac{\epsilon}{2e(H')}$ and $\beta = \frac{\epsilon^4}{e(H')}$. Let $L$ be large enough to invoke \cref{lemma: increment} on $\mathcal{H}(H')$. We obtain a set $W_0$ such that $|W_0| \leq \frac{\epsilon}{2e(H')}\binom{n}{2}$ and $\mathcal{E}_{H'}(W_0) \leq \frac{\epsilon^4}{e(H')}\binom{n}{v(H')}$ where $\mathcal{E}_{H'}$ is defined as in (\ref{eq: energy}). Now let $W \supseteq W_0$ be any set with $\leq \frac{\epsilon}{e(H')}\binom{n}{2}$ edges. Note that $\mathcal{E}_{H'}$ is a decreasing function with respect to containment, so $\mathcal{E}_{H'}(W) \leq \frac{\epsilon^4}{e(H')} \binom{n}{v(H')}$ as well. Applying \cref{lemma: cauchy schwarz}, we obtain 
    \begin{align*}
        \E{\sum_{A \in \mathcal{H}(H')-W} \abs{\E{Y_A\vert (Y_w)_{w \in W}} - \prod_{a \in A}\E{Y_a \vert (Y_w)_{w \in W}}}} &\leq p^{e(H')}(e(H')e(\mathcal{H}(H')-W))^{1/2}\mathcal{E}_{H'}(W)^{1/2} \\ &\leq \epsilon^2p^{e(H')}\binom{n}{v(H')}.
    \end{align*}
    By Markov's inequality, we have that with probability at least $1-\epsilon$, 
    \[ \sum_{A \in \mathcal{H}(H')-W}\prod_{a \in A}\E{Y_a \vert (Y_w)_{w \in W}} \leq \sum_{A \in \mathcal{H}(H')-W} \E{Y_A\vert (Y_w)_{w \in W}} + \epsilon p^{e(H')}\binom{n}{v(H')}. \]
    Since $Y$ is drawn from the lower tail distribution, it deterministically satisfies 
    \[ \sum_{A \in \mathcal{H}(H')-W} \E{Y_A\vert (Y_w)_{w \in W}} \leq \eta p^{e(H')}\binom{n}{v(H')} \]
    and so with probability $1-\epsilon$, we have 
    \[ \sum_{A \in \mathcal{H}(H')-W}\prod_{a \in A}\E{Y_a \vert (Y_w)_{w \in W}} \leq (\eta + \epsilon)p^{e(H')}\binom{n}{v(H')}. \]
    Define $q^W_e \vcentcolon= \E{Y_e \vert (Y_w)_{w \in W}}$ for $e \not \in W$ and $q^W_e = p$ for $e \in W$. The size of $W$ along with \cref{container condition} applied to $\mathcal{H}(H')$ implies that there are at most $\epsilon n^{v(H')}$ edges $A \in \mathcal{H}(H')$ such that $A \cap W \neq \emptyset$. \cref{lemma: harris} implies that $q_e^W \leq p$ for every vertex $e$. These two estimates combined ensure that 
    \[ \sum_{\substack{A \in \mathcal{H}(H') \\ A \cap W \neq \emptyset}}\prod_{a \in A}q^W_a \leq \epsilon p^{e(H')}\binom{n}{v(H')} \]
    and so all together we have 
    \[ \sum_{A \in \mathcal{H}(H')}\prod_{a \in A}q^W_a \leq (\eta + 2\epsilon)p^{e(H')}\binom{n}{v(H')}. \]
    In particular, $\mathbf{q}^W$ satisfies $\LT_p(H', \eta+2\epsilon)$ with probability at least $1-\epsilon$. 

    We use this fact to argue that $\q^W$ must be close to the constant vector $q$, where $q = \eta^{1/e(H')} p$. By the above we know that with probability at least $1-\epsilon$, 
    \[ i_p(\q^W) \geq \Phi_p(H', \eta + 2\epsilon). \] 
    Let $a(\delta) \vcentcolon= \Prob{i_p(\q^W) > \Phi_p(H', \eta) + \delta \binom{n}{2}p}$. This gives a lower bound of
    \[ a(\delta)(\Phi_p(H', \eta) + \delta \binom{n}{2}p) + (1-\epsilon - a(\delta))\Phi_p(H', \eta + 2\epsilon) \leq \E{i_p(\q^W)}. \]
    From Theorem \ref{theorem: tail probability} we can deduce an upper bound on $\E{i_p(\q^W)}$ via relative entropy manipulations. Using the properties from Definition \ref{def: entropy}, 
    \begin{align*}
        \E{i_p(\q^W)} &= \sum_{e \in \mathcal{V}\setminus W} I_p(Y_e | (Y_w)_{w \in W}) \leq I_p((Y_e)_{e \in \mathcal{V} \setminus W}|(Y_w)_{w \in W}) + I_p((Y_w)_{w \in W}) \\
        &\leq I_p(Y) = -\log \Prob{\LT_p(H', \eta)} \leq (1+\epsilon)\Phi_p(H', \eta(1-\epsilon)).
    \end{align*}
    Combining the two inequalities and rearranging for $a(\delta)$, 
    \[ a(\delta) \leq \frac{(1+\epsilon)\Phi_p(H', \eta(1-\epsilon)) - (1-\epsilon)\Phi_p(H', \eta + 2\epsilon)}{\Phi_p(H', \eta) + \delta \binom{n}{2}p - \Phi_p(H', \eta + 2\epsilon)} \lesssim \frac{\epsilon}{\epsilon + \delta}. \]
    Here we used the fact that each of the variational problems in the expression is minimized by the constant graphon, and the optimal value is locally Lipschitz in $\eta$. By \cref{prop: stability}, 
    \[ \Prob{\norm{\q^W - q}_\square > \delta p} \leq \alpha\left(\left(\frac{\delta}{C_{\eta, H}}\right)^{12}\right) \lesssim \frac{\epsilon}{\epsilon + \delta^{12}}. \] 
    Setting $\delta = \epsilon^{1/13}$ finishes the proof. 
\end{proof}
To boost this weak structure to a typical structure in terms of subgraph counts, we first estimate the expected count of any other graph $H$. This is done by once again invoking the entropy increment to gain control over copies of $H$ under the conditional distribution.
\begin{theorem}\label{thm: expectation main theorem}
    Fix a graph $H$ and let $m = \max\{m_2(H), m_2(H')\}$. Then for any $\epsilon > 0$ there exists an $L$ such that if $p \geq Ln^{-1/m}$, 
    \[ \abs{\E{N_H(G_\LT)} - \E{N_H(q)}} < \epsilon \E{N_H(q)}. \]
\end{theorem}
\begin{proof}
    Once again, we may assume that $\epsilon$ is sufficiently small as a function of $\eta, H'$ and $H$. Let $L$ be large enough so that the conditions of Lemma \ref{lemma: increment} are satisfied for both $H$ and $H'$ and the conditions of Lemma \ref{lemma: weak structure} are satisfied for $H'$. Let $W_0$ be the set guaranteed by \cref{lemma: weak structure} with parameter $\frac{\epsilon}{2e(H)e(H')}$. Invoking \cref{lemma: increment} on $\mathcal{H}(H)$ with $\alpha = \frac{\epsilon}{2e(H)e(H')}$ and $\beta = \frac{\epsilon^4}{e(H)}$, we obtain a set $W_1$ such that $|W_1| \leq \frac{\epsilon}{2e(H)e(H')} \binom{n}{2}$ and $\mathcal{E}_H(W_1) \leq \frac{\epsilon^4}{e(H)} \binom{n}{v(H)}$. Let $W = W_0 \cup W_1$, then $|W| \leq \frac{\epsilon}{e(H)e(H')}$ and $\mathcal{E}_H(W) \leq \frac{\epsilon^4}{e(H)}\binom{n}{v(H)}$. Moreover, by choice of $W_0$ and Lemma \ref{lemma: weak structure} we also have that 
    \[ \Prob{\norm{\q^W - q}_\square > \epsilon p} \leq C_{\eta, H'}\epsilon. \]

    Note that since $\q^W \leq p$ pointwise by \cref{lemma: harris}, the graphons $\frac{\q^W}{p}$ and $\frac{q}{p}$ are bounded above by 1. Thus, the counting lemma \cite[Section 4.5]{Zha:23} implies that 
    \[ \abs{\E{N_H(\q^W)} - \E{N_H(q)}} = p^{e(H)}\abs{\E{N_H\left(\frac{\q^W}{p}\right)} - \E{N_H\left(
    \frac{q}{p}\right)}} \leq p^{e(H)}n^{v(H)}e(H)\norm{\frac{\q^W}{p} - \frac{q}{p}}_\square. \]
    Here we treat $\q^W$ as fixed and the expectation is taken over the randomness of $N_H(G)$ where $G$ is a random graph sampled from $G(n, \q^W)$. This implies that 
    \[ \Prob{\abs{\E{N_H(\q^W)} - \E{N_H(q)}} > \epsilon e(H) p^{e(H)}n^{v(H)}} \leq \Prob{\norm{\q^W - q}_\square > \epsilon p} \leq C_{\eta, H'} \cdot \epsilon. \]

    Now we control the expected number of copies of $H$. We decompose
    \begin{align*}
        \E{N_H(Y)} &= \E{\sum_{A \in \mathcal{H}(H)} Y_A} = \E{\sum_{A \in \mathcal{H}(H)-W} Y_A} + \E{\sum_{A \in \mathcal{H}(H), A \cap W \neq \emptyset} Y_A}
    \end{align*}
    For the second expectation, note that the size of $W$ along with \cref{container condition} applied to $\mathcal{H}(H)$ implies that there are at most $\epsilon n^{v(H)}$ terms in the sum. Moreover, \cref{lemma: harris} implies that $\E{Y_A} \leq p^{e(H)}$ for every $A$, and thus the contribution of this term is bounded by $\epsilon p^{e(H)}n^{v(H)}$. 
    
    For the first expectation, applying \cref{lemma: cauchy schwarz} and an analogous computation to the proof of \cref{lemma: weak structure} to the hypergraph $\mathcal{H}(H)$, we have that with probability at least $1-\epsilon$,
    \[ \sum_{A \in \mathcal{H}(H)-W} \abs{\E{Y_A\vert (Y_w)_{w \in W}} - \prod_{a \in A}q_a^W} \leq \epsilon p^{e(H)}\binom{n}{v(H)}.  \]
    Moreover, recall that with probability at least $1-C_{\eta, H'} \epsilon$ we have 
    \[ \abs{\sum_{A \in \mathcal{H}(H)} \left(\prod_{a \in A}q_a^W - q^{e(H)}\right)} < \epsilon e(H) p^{e(H)}n^{v(H)}. \]
    Once again, since $\mathbf{q}^W$ and $q$ are both bounded by $p$, the size of $W$ implies that on this event we have 
    \[ \abs{\sum_{A \in \mathcal{H}(H)-W} \left(\prod_{a \in A}q_a^W - q^{e(H)}\right)} < 2\epsilon e(H) p^{e(H)}n^{v(H)} \]
    as well.
    Thus, with probability at least $1 - C_{\eta, H'}\epsilon$, 
    \[ \abs{\sum_{A \in \mathcal{H}(H)-W} \left(\E{Y_A\vert (Y_w)_{w \in W}} - q^{e(H)}\right)} \leq 3\epsilon e(H)p^{e(H)}n^{v(H)}. \]
    Note that by \cref{lemma: harris} we know that the random variable $\E{\sum_{A \in \mathcal{H}(H)-W} Y_A\middle\vert (Y_w)_{w \in W}}$ is uniformly bounded above by $|\mathcal{H}(H)-W|p^{e(H)} \leq n^{v(H)}p^{e(H)}$. Thus, 
    \begin{align*}
        \abs{\E{\sum_{A \in \mathcal{H}(H)-W} Y_A} - \binom{n}{v(H)}q^{e(H)}} &= \abs{\E{\E{\sum_{A \in \mathcal{H}(H)-W} Y_A\middle\vert (Y_w)_{w \in W}} - \binom{n}{v(H)}q^{e(H)}}} \\
        &\lesssim_{\eta, H, H'} \epsilon p^{e(H)}n^{v(H)}.
    \end{align*}
    Combining the estimates together, we obtain that
    \[ \abs{\E{N_H(Y)} - \binom{n}{v(H)}q^{e(H)}} \lesssim_{\eta, H,H'} \epsilon p^{e(H)}n^{v(H)} \]
    which is the desired conclusion.
\end{proof}
Finally, to prove \cref{thm: main intro}, we further boost the estimate on the expected number of copies of $H$ to a concentration of the number of copies of $H$ by a second moment argument. The key observation is that the variance can be written in terms of expected subgraph counts, which are controlled by \cref{thm: expectation main theorem}.

\begin{proof}[Proof of Theorem \ref{thm: main intro}]
    We show that the subgraph count $H$ concentrates well around its expectation, which by \cref{thm: expectation main theorem} implies that it is close to $\E{N_H(q)}$ as well. 
    \begin{align*}
        \Var{N_H(Y)} &= \Var{\sum_{A \in \mathcal{H}(H)} Y_A} = \sum_{A ,A' \in \mathcal{H}(H)} \E{Y_A Y_{A'}} - \E{Y_A}\E{Y_A'}
    \end{align*}
    We partition this sum into cases based on the graph $\Gamma$ induced by the edge set $A \cup A'$. Suppose first that $A \cap A' \neq \emptyset$. Denote the graph of their intersection by $\Gamma'$, and note that $\Gamma'$ is a subgraph of $H$. By \cref{lemma: harris} we can bound
    \begin{align*}
        \E{Y_AY_{A'}} \leq p^{e(\Gamma)} = p^{2e(H) - |e(\Gamma')|}
    \end{align*}
    Moreover, there are at most $n^{v(\Gamma)} = n^{2v(H) - v(\Gamma')}$ pairs of $A$ and $A'$ that induce this graph. Thus, the total contribution is bounded by 
    \[ \frac{n^{2v(H)}p^{2e(H)}}{n^{v(\Gamma')}p^{e(\Gamma')}} \ll n^{2v(H)}p^{2e(H)} \]
    for the given regime of $p$. Thus, it only remains to consider $A \cap A' = \emptyset$, in other words $\Gamma =$ the disjoint union of two copies of $H$. Note that $m_2(\Gamma) = m_2(H)$ and thus by \cref{thm: expectation main theorem} we know that for $L$ large enough $\abs{\E{N_\Gamma(Y)} - \E{N_\Gamma(q)}} < \epsilon^3 p^{e(\Gamma)}n^{v(\Gamma)}$. By the symmetry of the copies of $\Gamma$, this implies that $\abs{\E{Y_A Y_{A'}} - q^{e(\Gamma)}} < \epsilon^3 p^{e(\Gamma)}$. In particular, for the pairs $A \cap A' = \emptyset$ the contribution is bounded by $\epsilon^3 p^{e(\Gamma)}$. There are at most $n^{2v(H)}$ pairs in the sum, so
    \[ \Var{N_H(Y)} \leq \epsilon^3 p^{2e(H)}n^{2v(H)} \leq \epsilon^3 \E{N_H(Y)}^2. \]
    The result follows from Chebyshev's inequality. 
\end{proof}

To deduce Corollary \ref{cor: main cor} we show that there exist finite graphs that certify pseudorandomness with 2-density arbitrarily close to 1. As a consequence typicality in the sense of cut norm follows from typicality in the sense of subgraph counts.

\begin{proof}[Proof of Corollary \ref{cor: main cor}]
    Choose $k$ large enough such that $\frac{k}{k-1} < m_2(H)$. Since $m_2(C_{2k}) = \frac{2k-1}{2k-2} \leq \frac{k}{k-1} < m_2(H)$, we may apply \cref{thm: main intro} to $H$ and $C_{2k}$ and obtain some constant $L$ such that for $p \geq Ln^{-1/m_2(H)}$, $G_\LT$ has the correct $C_{2k}$ count with probability at least $1-\epsilon$. We show that $C_{2k}$ is forcing for all $p \geq Ln^{-\frac{1}{m_2(H)}} \gg n^{-\frac{k-1}{k}}$, in particular  the correct $C_{2k}$ count implies pseudorandomness (see \cite[Chapter 3]{Zha:23} for further discussion of forcing graphs). This implies that $G_\LT$ is close to the constant graphon with probability at least $1-\epsilon$.

    Let $A$ be the adjacency matrix of $G_\LT$. We know that with probability at least $1-\epsilon$, 
    \[ \abs{N_{C_{2k}}(G_\LT) - \binom{n}{2k}q^{2k}} < \epsilon n^{2k}p^{2k}. \] 
    Note that $\Tr(A^{2k})$ counts the number of closed walks of length $2k$ in the graph $G_\LT$. Let $H$ be any graph that can be a closed walk of length $2k$. By Harris's inequality $\E{N_H(G_\LT)} \leq \binom{n}{v(H)}p^{e(H)}$. This is much smaller than $n^{2k}p^{2k}$ as soon as $p \gg n^{-\frac{2k-v(H)}{2k-e(H)}}$. This quantity is maximized when $H$ is a tree with $k$ edges (each of which is traversed twice to form a closed walk). Since $p \gg n^{-\frac{k-1}{k}}$, for any such graph $H$ the expected contribution is $o(n^{2k}p^{2k})$, and thus with probability $1-o(1)$ $N_H(G_\LT) = o(n^{2k}p^{2k})$. It remains to show that $\norm{A - qJ}_\square < \epsilon p$ where $J$ is the all 1's matrix. The above computation verifies that $\abs{\Tr((A-qJ)^{2k})} < 2\epsilon n^{2k}p^{2k}$ with probability $1-2\epsilon$ and on this event the eigenvalues of $A - qJ$ are all at most $\epsilon^{1/2k}np$. For any vectors $x, y \in [0,1]^n$, write $x = \sum_{i=1}^n a_i v_i$ and $y = \sum_{i=1}^n b_i v_i$ where $\{v_i\}$ is an orthonormal basis of eigenvectors of $A - qJ$. Then 
    \begin{align*}
        \abs{x^\top (A - qJ)y} &= \abs{\sum_{i=1}^n a_ib_i\lambda_i} \leq (\max_i |\lambda_i|)\sum_{i=1}^n |a_ib_i| \leq o(np)\cdot \left(\sum_{i=1}^n a_i^2\right)^{1/2}\left(\sum_{i=1}^n b_i^2\right)^{1/2} \\
        &= \epsilon^{1/2k}np \cdot \norm{x}_2\norm{y}_2 \leq \epsilon^{1/2k}n^2p.
    \end{align*}
    The conclusion follows since $\norm{G_\LT - q}_\square = \sup_{x,y \in [0,1]^n} \frac{1}{n^2}\abs{x^\top (A - qJ)y} = \epsilon^{1/2k}p$ and $\epsilon > 0$ was chosen arbitrarily. 
\end{proof}

\section{Discussion: extension to general hypergraphs}\label{section: discussion}
To conclude, we discuss the applicability of this argument in the more general hypergraph setup of \cite{KS:23}. We first introduce the notation that extends our study of lower tails in hypergraphs encoding subgraphs -- namely $\mathcal{H}(H)$ -- to general hypergraphs. Recall that we denote a hypergraph by $\mathcal{H} = (\mathcal{V}, \mathcal{E})$. First, we define the corresponding lower tail event, which says that a random subset of the vertices contains fewer edges than expected.
\begin{definition}[General lower tail event]
    For $\mathbf{q} \in [0,1]^\mathcal{V}$ let $\mathcal{V}_\mathbf{q}$ be a random subset of $\mathcal{V}$ where each vertex $v$ is retained independently with probability $q_v$. When $q_v = p$ for every $v$ we write $\mathcal{V}_p$. The lower tail event $\LT_p(\mathcal{H}, \eta)$ is the event that $e(\mathcal{H}[U]) \leq \eta \E{e(\mathcal{H}[\mathcal{V}_p])}$. Here $\mathcal{H}[U]$ denotes the restriction of the hypergraph $\mathcal{H}$ to the vertex subset $U \subset \mathcal{V}$. 
\end{definition}
We can analogously define the variational problem associated with this lower tail event.

\begin{definition}[General variational problem]
    The variational problem associated with the lower tail event $\LT_p(\mathcal{H}, \eta)$ is
    \[ \Phi_p(\mathcal{H}, \eta) \vcentcolon= \min_{\mathbf{q} \in [0,1]^\mathcal{V}} \{ i_p(\mathbf{q}): \E{e(\mathcal{H}[\mathcal{V}_\mathbf{q}])} \leq \eta \E{e(\mathcal{H}[\mathcal{V}_p])}\} \]
\end{definition}

All of the results of \cite{KS:23} that we recounted in Section \ref{section: entropy increment} have analogs that hold for hypergraphs satisfying the following uniformity of degrees condition. In particular, Kozma and Samotij show that Theorem \ref{theorem: tail probability} holds for hypergraphs satisfying Condition \ref{container condition} -- that the logarithm of the lower tail probability is controlled by the associated variational problem.
\begin{condition}\label{container condition}
    For an $r$-uniform hypergraph $\mathcal{H}$, $v(\mathcal{H}) = |V(\mathcal{H}|$, $e(\mathcal{H}) 
    = |E(\mathcal{H})|$, and $\deg_\mathcal{H}(B) = |\{e \in E(\mathcal{H}): B \subset e\}|$. Let $\Delta_s \vcentcolon= \max\{\deg_\mathcal{H}(B): |B| = s\}$ be the maximum degree of a set of $s$ vertices. Then there exists $\lambda$ and $K$ such that the following holds for every $s \in [r]$: 
    \[ \Delta_s(\mathcal{H}) \leq K(\lambda p)^{s-1} \frac{e(\mathcal{H})}{v(\mathcal{H})}. \]
\end{condition}
The condition is closely related to that of the hypergraph container method, introduced independently by Balogh, Morris, and Samotij \cite{BMS:15} and Saxton and Thomason \cite{ST:15}, see \cite{BMS:18} for a survey of this powerful new technique. For further discussion of this connection, see the introduction of \cite{KS:23}. 

In this general setup, our argument extends in the following manner. Suppose we have two hypergraphs $\mathcal{H}$ and $\mathcal{H}'$ that live on the same vertex set and satisfy Condition \ref{container condition} under the same parameters.  We can apply the entropy increment to decompose the lower tail measure with respect to both $\mathcal{H}$ and $\mathcal{H}'$ simultaneously. If we know that $\Phi_p(\mathcal{H}', \eta)$ has a stability of the minimizing vector $\mathbf{q}$, we can deduce that the expected number of edges of $\mathcal{H}$ in a random sample of vertices from $\mathcal{V}_p$ conditioned on $\LT_p(\mathcal{H}', p)$ is close to the expectation under $\mathbf{q}$ (analogous to Theorem \ref{thm: expectation main theorem}). To boost this to a concentration statement we need to understand the pairwise interactions between edges $\mathcal{H}$. Under the further assumption that pairs of intersecting edges in $\mathcal{H}$ also satisfy the uniformity Condition \ref{container condition}, we obtain that the number of edges of $\mathcal{H}$ in a random sample of vertices from the conditional measure is close to the expectation under $\mathbf{q}$ with high probability (analogous to Theorem \ref{thm: main intro}). This reduces the problem of counting edges sampled from a lower tail measure to analyzing the variational problem and the uniformity of pairs of edges in the hypergraph.

The argument is also robust to conditioning on multiple lower tail events simultaneously, as well as counting edges of multiple hypergraphs simultaneously. However, given the difficulty of resolving the variational problem already for a single lower tail, applications of this form seem farther from fruition.

\bibliographystyle{amsplain0}
\bibliography{lowertails}

\providecommand{\bysame}{\leavevmode\hbox to3em{\hrulefill}\thinspace}
\providecommand{\MR}{\relax\ifhmode\unskip\space\fi MR }
\providecommand{\MRhref}[2]{%
  \href{http://www.ams.org/mathscinet-getitem?mr=#1}{#2}
}
\providecommand{\href}[2]{#2}
\begin{thebibliography}{10}

\bibitem{BMS:15}
J\'{o}zsef Balogh, Robert Morris, and Wojciech Samotij, \emph{Independent sets
  in hypergraphs}, J. Amer. Math. Soc. \textbf{28} (2015), 669--709.

\bibitem{BMS:18}
J\'{o}zsef Balogh, Robert Morris, and Wojciech Samotij, \emph{The method of
  hypergraph containers}, Proceedings of the {I}nternational {C}ongress of
  {M}athematicians---{R}io de {J}aneiro 2018. {V}ol. {IV}. {I}nvited lectures,
  World Sci. Publ., Hackensack, NJ, 2018, pp.~3059--3092.

\bibitem{BCCZ:19}
Christian Borgs, Jennifer~T. Chayes, Henry Cohn, and Yufei Zhao, \emph{An
  {$L^p$} theory of sparse graph convergence {I}: {L}imits, sparse random graph
  models, and power law distributions}, Trans. Amer. Math. Soc. \textbf{372}
  (2019), 3019--3062.

\bibitem{CD:16}
Sourav Chatterjee and Amir Dembo, \emph{Nonlinear large deviations}, Adv. Math.
  \textbf{299} (2016), 396--450.

\bibitem{CV:11}
Sourav Chatterjee and S.~R.~S. Varadhan, \emph{The large deviation principle
  for the {E}rd{\H{o}}s-{R}\'{e}nyi random graph}, European J. Combin.
  \textbf{32} (2011), 1000--1017.

\bibitem{CG:16}
D.~Conlon and W.~T. Gowers, \emph{Combinatorial theorems in sparse random
  sets}, Ann. of Math. (2) \textbf{184} (2016), 367--454.

\bibitem{CD:20}
Nicholas Cook and Amir Dembo, \emph{Large deviations of subgraph counts for
  sparse {E}rd{\H{o}}s-{R}\'{e}nyi graphs}, Adv. Math. \textbf{373} (2020),
  107289, 53.

\bibitem{CDP:21}
Nicholas~A. Cook, Amir Dembo, and Huy~Tuan Pham, \emph{Regularity method and
  large deviation principles for the erd{\H o}s--r{\'e}nyi hypergraph},
  (2021).

\bibitem{Eld:18}
Ronen Eldan, \emph{Gaussian-width gradient complexity, reverse log-{S}obolev
  inequalities and nonlinear large deviations}, Geom. Funct. Anal. \textbf{28}
  (2018), 1548--1596.

\bibitem{FK:99}
Alan Frieze and Ravi Kannan, \emph{Quick approximation to matrices and
  applications}, Combinatorica \textbf{19} (1999), 175--220.

\bibitem{JW:16}
Svante Janson and Lutz Warnke, \emph{The lower tail: {P}oisson approximation
  revisited}, Random Structures Algorithms \textbf{48} (2016), 219--246.

\bibitem{KRRS:18}
Richard Kenyon, Charles Radin, Kui Ren, and Lorenzo Sadun, \emph{Bipodal
  structure in oversaturated random graphs}, Int. Math. Res. Not. IMRN (2018),
  1009--1044.

\bibitem{KS:23}
Gady Kozma and Wojciech Samotij, \emph{Lower tails via relative entropy}, Ann.
  Probab. \textbf{51} (2023), 665--698.

\bibitem{Lov:12}
L\'{a}szl\'{o} Lov\'{a}sz, \emph{Large networks and graph limits}, American
  Mathematical Society Colloquium Publications, vol.~60, American Mathematical
  Society, Providence, RI, 2012.

\bibitem{NRS:23}
Joe Neeman, Charles Radin, and Lorenzo Sadun, \emph{Typical large graphs with
  given edge and triangle densities}, Probab. Theory Related Fields
  \textbf{186} (2023), 1167--1223.

\bibitem{Rad:18}
Charles Radin, \emph{Phases in large combinatorial systems}, Ann. Inst. Henri
  Poincar\'{e} D \textbf{5} (2018), 287--308.

\bibitem{ST:15}
David Saxton and Andrew Thomason, \emph{Hypergraph containers}, Invent. Math.
  \textbf{201} (2015), 925--992.

\bibitem{Sch:16}
Mathias Schacht, \emph{Extremal results for random discrete structures}, Ann.
  of Math. (2) \textbf{184} (2016), 333--365.

\bibitem{Zha:17}
Yufei Zhao, \emph{On the lower tail variational problem for random graphs},
  Combin. Probab. Comput. \textbf{26} (2017), 301--320.

\bibitem{Zha:23}
Yufei Zhao, \emph{Graph theory and additive combinatorics---exploring structure
  and randomness}, Cambridge University Press, Cambridge, 2023.

\end{thebibliography}

\end{document}